\documentclass[12pt]{article}
\usepackage{verbatim}
\usepackage{hyperref,url,breakurl}
\newcommand{\nocopyright}{
No Copyright\thanks{
The authors hereby waive all copyright
and related or neighboring rights to this work,
and dedicate it to the public domain.
This applies worldwide.
}}
\title{Cyclic groups with the same Hodge series}
\author{Daryl R. DeFord \and Peter G. Doyle}
\date{Version dated 8 April 2014
\\ \nocopyright}
\usepackage{amsthm,amssymb}

\newcommand{\divides}{\,|\,}

\newcommand{\LMR}{\mathrm{LMR}}

\newcommand{\diag}{\mathrm{diag}}
\newcommand{\dub}[1]{{#1^\pm}}
\newcommand{\Ldub}{L^\pm}
\newcommand{\quot}[1]{/_#1\,}
\newcommand{\eqmd}[1]{\equiv_{S_m \times\ZZ_{#1}}}

\newcommand{\bra}{\langle}
\newcommand{\ket}{\rangle}

\newcommand{\implies}{\Rightarrow}
\newcommand{\ZZstar}{\ZZ^\star}

\newcommand{\hodge}{\Lambda}

\newcommand{\intersect}{\cap}
\newcommand{\goesinto}{\backslash}
\newcommand{\goesto}{\rightarrow}
\newcommand{\into}{\rightarrow}
\newcommand{\CC}{\mathbf{C}}
\newcommand{\RR}{\mathbf{R}}
\newcommand{\ZZ}{\mathbf{Z}}
\newcommand{\eqhodge}{\equiv_{\hodge}}
\newtheorem{thm}{Theorem}

\newtheorem{lemma}{Lemma}
\newtheorem{prop}{Proposition}

\newtheorem{definition}{Definition}

\newtheorem{exercise}{Exercise}
\newcommand{\eqstart}{\begin{eqnarray*}}
\newcommand{\eqend}{\end{eqnarray*}}
\newcommand{\proofstart}{{\bf Proof.\ }}
\newcommand{\mathproofend}{\quad \qed}
\newcommand{\proofend}{$\quad \qed$}
\begin{document}
\maketitle
\begin{abstract}
The Hodge series of a finite matrix group 
is the generating function $\sum_{k,p} x^k y^p$ for invariant exterior forms
of specified order $p$ and degree $k$.
Lauret, Miatello, and Rossetti gave examples of pairs of
non-conjugate cyclic groups having the same Hodge series;
the corresponding space forms are isospectral for the
Laplacian on $p$-forms for all $p$, but not for all natural operators.
Here we explain, simplify, and extend their investigations.
\end{abstract}

\setcounter{section}{-1}
\section{Terminology and notation}

We adopt terminology
and notation to avoid some common headaches.

\subparagraph*{`Just if'.}
We follow John Conway in using `just if'
in place of the more cumbersome `if and only if'.

\subparagraph*{Modular arithmetic.}
$a \equiv_q b$ means $a$ is equivalent to $b$ mod $q$.
We write $\ZZ_q$  for $\ZZ/q\ZZ$, and $\ZZstar_q$ for its invertible elements,
taking $\ZZstar_1 = \{0\}$.
For $a \in \ZZ_q$, $b \in \ZZstar_q$
we write $a \quot{q} b$ for the quotient mod $q$.

\subparagraph*{Angles; roots of unity.}
We use $\tau=2\pi$ in representing angles,
because as Vi Hart
\cite{vihart:pi}
has so persuasively argued, \emph{$\pi$ is wrong}.
We write
\[
\omega_q = \exp(i \tau/q)
\]
for the standard $q$th root of unity, so that
\[
\omega_q^k = \exp(i \tau k/q) = e^{i \tau \frac{k}{q}}
.
\]

\subparagraph*{Unitary and orthogonal groups; conjugacy.}
As usual we write $U_n \subset GL_n(\CC)$ for the 
$n$-by-$n$ unitary matrices, and $O_n = U_n \intersect GL_n(\RR)$ for
the orthogonal matrices.

When we say that two matrices or groups are `conjugate', we
mean that they are conjugate within $GL_n(\CC)$, so that the
conjugating matrix can be any invertible complex matrix.
Allowing this generality for the conjugating matrix is no big deal, because
unitary matrices or groups that are conjugate within $GL_n(\CC)$ are already
conjugate within $U_n$;
real matrices or groups that are conjugate within $GL_n(\RR)$ are already
conjugate within $GL_n(\RR)$;
orthogonal matrices or groups that are conjugate within $GL_n(\CC)$ are already
conjugate within $O_n$.

We will be dealing with finite groups of matrices, which we
will be interested in only up to conjugacy.
Any finite group of complex matrices is conjugate to a subgroup of $U_n$;
any finite group of real matrices is conjugate to a subgroup of $O_n$.
So we may take our groups to be unitary---and if real, orthogonal---without
sacrificing generality.

\section{Hodge series}

Let $G \subset U_n$ be a finite group of $n$-by-$n$ complex matrices,
assumed to be unitary.
Any $g \in G$ is diagonalizable,
with the roots $\lambda_1,\ldots,\lambda_n$ of its characteristic polynomial
\[
\chi_g(x) = \det(I_n-xg) = \prod_i (x-\lambda_i)
\]
being roots of unity.

Define the \emph{Hodge series}
\eqstart
\hodge_G(x,y)
&=&
\frac{1}{|G|}\sum_g \frac{\det(I_n+yg)}{\det(I_n-xg)}
\\&=&
\frac{1}{|G|}\sum_g \frac{y^n \chi_g(-1/y)}{x^n \chi_g(1/x)}
.
\eqend
This series is a particular kind of \emph{Molien series}:
Crass \cite[p.\ 31]{crass:sextic} calls it the `exterior Molien series'.
By a generalization of Molien's theorem
(cf. Molien \cite{molien}, Stanley \cite{stanley:molien})
this is the generating function for
$G$-invariant exterior forms:
\[
	\hodge_G(x,y) = \sum_{p,k} P_k^p x^k y^p
	,
\]
where $P_k^p$ is the dimension of the space of $G$-invariant $p$-forms
whose coefficients are homogeneous polynomials of degree $k$
in $x_1,\ldots,x_n$.

For example we have
\[
	\hodge_{\{I_n\}} = \frac{(1+y)^n}{(1-x)^n}
\]
and
\eqstart
\hodge_{\{\pm I_n\}}
&=&
\frac12
\left(\frac{(1+y)^n}{(1-x)^n}+\frac{(1-y)^n}{(1+x)^n}\right)
\\&=&
\frac{\frac12 ((1+x^n)(1+y^n)+(1-x^n)(1-y^n))}{(1-x^2)^n}
.
\eqend

{\bf Aside.}
Here's a more interesting example.
The group $G_{120}$ of
proper and improper symmetries
of the icosahedron in Euclidean 3-space has Hodge series
\[
\hodge_{G_{120}} =
\frac
{ (1 + x y) (1 + x^5 y) (1 + x^9 y)}
{(1 - x^2) (1 - x^6) (1 - x^{10})}
.
\]
As this might suggest,
the algebra of invariant forms is generated by polynomial
invariants of degrees $2,6,10$ and their exterior derivatives
of degrees $1,5,9$.
For the index-2 subgroup $G_{60}$ of proper symmetries we have
\[
\hodge_{G_{60}} = 
\frac
{ (1 + x^{15})(1+y^3) + (x + x^5 + x^6 + x^9 + x^{10} + x^{14}) (y+y^2)}
{(1 - x^2) (1 - x^6) (1 - x^{10})}
.
\]
This is a little harder to decipher, though the generating function
for invariant polynomials, obtained by setting $y=0$, is clear enough:
\[
\hodge_{G_{60}}(x,0) =
\frac
{ 1 + x^{15}}
{(1 - x^2) (1 - x^6) (1 - x^{10})}
.
\]
Here we see the $G_{120}$-invariants of degrees 2,6,10,
together with a new invariant of degree $15$ (the product of the linear
forms determining the $15$ planes of symmetry of the icosahedron)
whose square is $G_{120}$-invariant,
though it itself is only $G_{60}$-invariant.

\begin{exercise}
Compute these two Hodge series.
\end{exercise}

{\bf Hint.}
Resist the temptation to consult
Klein
\cite{klein:ikos}
or
Doyle and McMullen
\cite{doylemcmullen:icos}:
You do not need to know the matrix groups explicitly, because
the contribution of a matrix to the Hodge series depends only
on its conjugacy class.
This fact is the basis of the notion of `almost-conjugacy' of
groups, which we'll get to in a jiffy.

\section{Hodge equivalence}

We are interested in pairs of groups $G, H \subset U_n$
(and in particular, pairs of real groups $G,H \subset O_n$)
having the same Hodge series,
meaning that they have the same dimensions of spaces of invariant forms.
We call such pairs \emph{Hodge-equivalent}, and write
$G \eqhodge H$.

Of course if the groups $G$ and $H$ are conjugate
then they are Hodge-equivalent.
More generally,
say that $G$ and $H$ are \emph{almost conjugate} if there is a bijection
$\sigma:G \into H$ such that $g$ and $\sigma(g)$ are conjugate.
This is the same as requiring that $g$ and $\sigma(g)$
have the same eigenvalues,
so that they make the same
contribution to the Hodge series.
Thus almost conjugate groups are Hodge-equivalent.

It is relatively easy to find non-conjugate
pairs $(G,H)$ that are almost conjugate,
and hence Hodge-equivalent.
(See Gilkey \cite{gilkey:metacyclic},
Ikeda
\cite{ikeda:subtle}.)
Here we are interested in pairs (and specifically, real pairs)
that are Hodge-equivalent without being almost conjugate.
The first such examples were given 
in
\cite{lmr:lens}
by
Lauret, Miatello, and Rossetti 
(henceforth `LMR').
They exhibited a multitude of examples arising already among cyclic groups.
For cyclic groups, conjugacy is the same as almost-conjugacy,
so their examples can be briefly described as being Hodge-equivalent without
being conjugate.
Our goal here is to explain, simplify, and extend their findings.

\section{Isospectrality}

We discuss here the connection to spectral theory, which is what
motivated LMR to construct their examples.
This is meant for background only:
In the approach taken here, spectral theory plays no role.
In this section we restrict to real groups, which we may 
assume to be orthogonal.

A finite real group $G \subset O_n$ is classified up to conjugacy
by the isometry type of the quotient orbifold $Q_G = G \goesinto S^{n-1}$.
According to Ikeda
\cite{ikeda:subtle},
$G$ and $H$ are Hodge-equivalent
just if
the quotients
$Q_G$ and $Q_H$ are isospectral for the Hodge Laplacian on $p$-forms for
$p=0,\ldots,n-1$.
According to Pesce
\cite{pesce:strong},
$Q_G$ and $Q_H$ are \emph{strongly isospectral} 
(isospectral for all natural operators of a certain kind)
just if $G$ and $H$ are almost conjugate.
Using this dictionary, looking for Hodge-equivalent groups
that are not almost conjugate is the same as looking for Hodge-isospectral
orbifolds that are not strongly isospectral.
This is why LMR were interested in this question.

Ikeda, Pesce, and LMR restricted their investigations
to the case of groups whose action on $S_{n-1}$ is fixed-point free
(no $g \neq 1$ has $1$ as an eigenvalue).
In this case $Q_G$ is a manifold, called
a \emph{spherical space form}.
For $n$ odd (i.e. $n-1$ even) we have only
the sphere $\{I_n\} \goesinto S^{n-1}$
and projective space $\{\pm I_n \} \goesinto S^{n-1}$.
So the restriction to fixed-point free actions effectively limits us to
the case of even $n$.

When $G$ is cyclic, as in the LMR examples,
$Q_G$ is a \emph{lens space}.
As observed above, for cyclic groups
almost-conjugacy is the same as conjugacy,
which is the same as isometry of the corresponding lens space.
So we can briefly describe the LMR examples as 
Hodge-isospectral lens spaces that are not isometric,
hence not almost conjugate, hence not strongly isospectral.

This ends our discussion of isospectrality.
The rest is algebra.

\section{Cyclic groups}

For any $q$ and $s =(s_1,\ldots,s_n)$,
write
\[
\omega_q^s = (\omega_q^{s_1},\ldots,\omega_q^{s_n})
.
\]
Consider the finite cyclic group
\[
L(q,s)
=
\bra \diag(\omega_q^s) \ket = \{ \diag(\omega_q^{ks}): k \in \ZZ_q \}
.
\]
This group has order $q$ if $\gcd(q,s_1,\ldots,s_n) = 1$.
Up to conjugacy,
any finite cyclic subgroup of $U_n$ can be written in this way.

The cyclic group $L(q,s)$ doesn't change (up to conjugacy)
when you rearrange the entries of $s$, or multiply them
all by an element of the multiplicative group $\ZZstar_q$.
Conversely, the groups $L(q,s)$ and $L(q,s')$ are conjugate
just if,
when viewed as multisets mod $q$,
$s'$ can be obtained from $s$ by 
multiplying by an invertible element.

Now take $n=2m$, and let $\rho:GL_m(\CC) \mapsto GL_{2m}(\RR)$ be
the standard embedding,
so that
$\rho(\omega_q^s)$ is the diagonal sum of the 2-by-2 matrices
\begin{eqnarray*}
\rho(((\omega_q^{s_i})))
&=&
\exp(\tau \frac{s_i}{q} ((0,-1),(1,0)))
\\&=&
((\cos(\tau s_i/q),-\sin(\tau s_i/q)),(\sin(\tau s_i/q),\cos(\tau s_i/q)))
.
\end{eqnarray*}

Up to conjugacy in $U_n$,
\[
\rho(L(q,s))
\equiv L(q,\dub{s})
\]
where
\[
\dub{s} = 
(s_1,-s_1,\ldots,s_m,-s_m)
.
\]
Let us write
\[
\Ldub(q,s) = L(q,\dub{s}) \equiv \rho(L(q,s))
.
\]
For the Hodge series we have
\eqstart
\hodge_{\Ldub(q,s)}
&=&
\frac{1}{q}
\sum_k \prod_i
\frac{(1+y \omega_q^{k s_i})(1+y \omega_q^{-k s_i})}
{(1-x \omega_q^{k s_i})(1-x \omega_q^{-k s_i})}
\\&=&
\frac{1}{q}
\sum_k \prod_i
\frac{1+ 2 \cos(\tau k s_i/q) y + y^2}{1 - 2 \cos(\tau k s_i/q) x +x^2}
.
\eqend

\section{The LMR construction}

The LMR examples involve cyclic subgroups of $O_{2m}$ of the form
\[
\rho(L(r^2t,rta+1)) \equiv \Ldub(r^2t,rta+1)
,
\]
where
$r>2$, $t \geq 1$, $a=(a_1,\ldots,a_m) \in \ZZ^m$.
Since we prefer to keep our matrices diagonal we'll define
\eqstart
\LMR(r,t,a)
&=&
\Ldub(r^2t,rta+1)
\\&=&
L(r^2t,(rt a_1+1,-rt a_1 -1,\ldots,rt a_m +1,-rt a_m -1))
.
\eqend

{\bf Note.}  You may wish to mentally set $t=1$:
All evidence indicates that what works for $t=1$ works in general,
and in particular the criterion in Theorem \ref{one} below
does not involve $t$.

As we will be seeing, what's special about the LMR construction is the
following fact:
\[
(rtc+1)(rtd+1) \equiv_{r^2t} rt(c+d)+1
.
\]
Thus the multiplicative subgroup
$\{rtc+1:c \in \ZZ_r\} \subset \ZZstar_{r^2t}$
is cyclic of order $r$,
generated by $rt+1$;
the map $rtc+1 \mapsto c$
takes the logarithm of $rtc+1$ base $rt+1$,
and gives an isomorphism to the additive group $\ZZ_r$.

As a first consequence of this,
notice that we can add a constant $c$ to the entries
of $a$ without changing the conjugacy class:
\[
\LMR(r,t,a) \equiv \LMR(r,t,a+c)
.
\]
In fact, this characterizes all such coincidences:
Let us write
\[
a \eqmd{r} a'
\]
if for some $c$, $a+c$ and $a'$ are the same as multisets mod $r$.
Then $\LMR(r,t,a) \equiv \LMR(r,t,a')$
just if
$a \eqmd{r} a'$.

All the LMR pairs are (conjugate to) pairs of the special form
\[
(\LMR(r,t,a), \LMR(r,t,-a))
.
\]
Not all such pairs are Hodge-equivalent, however.

\section{Theorem}

In this section
we formulate a criterion for Hodge-equivalence of the LMR pair
$(\LMR(r,t,a), \LMR(r,t,-a))$.
While this criterion has not been shown to be necessary,
it holds in all the cases (thousands and thousands!)
where the LMR construction has been found to succeed.

\begin{definition}
Say that
$a=(a_1,\ldots,a_m)$ is:
\begin{itemize}
\item
\emph{univalent mod $r$} if its entries are distinct mod $r$;
\item
\emph{reversible mod $r$} if
$a \eqmd{r} -a$;
\item
\emph{good mod $r$} if it is univalent or reversible mod $r$;
\item
\emph{hereditarily good mod $r$} if it is good mod $d$ for all $d$ dividing $r$;
\item
\emph{useful mod $r$} if it is hereditarily good and irreversible mod $r$.
\end{itemize}
\end{definition}

Any $a$ is reversible (hence good) mod $1$ or $2$.
So in checking hereditary goodness we need only check divisors $d>2$.

In section \ref{pf} below we will prove the following:
\begin{thm} \label{one}
If $a$ is hereditarily good mod $r$ then for any $t$,
\[
\LMR(r,t,a) \eqhodge \LMR(r,t,-a)
.
\]
\end{thm}

If $a$ is reversible mod $r$ then $a$ is hereditarily good mod $r$,
but in this case $\LMR(r,t,a)$ and $\LMR(r,t,-a)$ are conjugate.
So this result tells us something useful only if $a$ is hereditarily good
without being reversible, which is our definition of `useful'.

\section{Examples}
$(0,1,3)$ is:
\begin{itemize}
\item
univalent mod $4,5,6,\ldots$;
\item
reversible mod $1,2,4,5$;
\item
good mod any $r \neq 3$;
\item
hereditarily good mod any $r$ not divisible by $3$;
\item
useful mod any $r \geq 7$ not divisible by $3$.
\end{itemize}
Putting $r=7,8,10$, $t=1$, we get 
Hodge-equivalent but non-conjugate pairs of orders $49,64,100$;
Putting $r=7$, $t=2$ we get a pair of order $98$.

\medskip
$(0,1,4)$ is:
\begin{itemize}
\item
univalent mod $5,6,7,\ldots$;
\item
reversible mod $1,2,5,7$;
\item
good mod any $r \neq 3,4$;
\item
hereditarily good mod any $r$ not divisible by $3$ or $4$;
\item
useful mod any $r \geq 10$ not divisible by $3$ or $4$.
\end{itemize}
Putting $r=10$, $t=1$ gives a pair of order 100.
Together with the four pairs coming from $(0,1,3)$ above,
this gives us all five inequivalent pairs with $m=3$, $q \leq 100$
(see Table 1 of LMR \cite{lmr:lens}).

We'll call the simplest of these pairs the \emph{$49$-pair}:
\begin{eqnarray*}
&&
(\LMR(7,1,(0,1,3)),\LMR(7,1,(0,-1,-3)))
\\&=&
(\Ldub(49,(1,8,22)),\Ldub(49,(1,-6,-20))
\\&=&
(L(49,(1,-1,8,-8,22,-22)),L(49,(1,-1,-6,6,-20,20)))
\\&\equiv&
(L(49,(-6,6,1,-1,15,-15)),L(49,(1,-1,-6,6,-20,20)))
\\&\equiv&
(\Ldub(49,(1,6,15)),\Ldub(49,(1,6,20)))
.
\end{eqnarray*}
Here at the next-to-last step we've multiplied the list
$(1,-1,8,-8,22,-22)$ by $-6$ mod $49$ so as to get the lexicographically
least representation that the computer spits out in its search
for Hodge-equivalent pairs.

\section{Proof} \label{pf}

It is easy enough to verify that the members of the $49$-pair
are Hodge-equivalent
by explicit computation of their Hodge series.
The same goes for as many other pairs as you like, but this only gets you
a finite number of examples.

Using a very explicit representation theory argument,
LMR proved Hodge-equivalence of the 49-pair in a way that extends
to cover all pairs of the form
\[
(\LMR(r,t,(0,1,3)),\LMR(r,t,(0,-1,-3)))
\]
with $r$ not divisible by $3$.
As we have seen, this
infinite family
is just what we get out of Theorem \ref{one} if we take $a=(0,1,3)$.
It includes $19$ of the $62$ examples
in the list given by LMR of all pairs with $m=3$ and $q \leq 300$.

To prove Theorem \ref{one} in its full generality,
we're going to show that the two Hodge series involved
are identical as rational
functions of $x$ and $y$.
This comes down to a bunch of manipulations with partial fraction expansions.
It all starts with the following familiar identity.

\begin{lemma}
\[
\prod_{i=1}^n \frac{1}{x-\lambda_i} =
\sum_{i=1}^{n}
\frac{1}{x-\lambda_i} \prod_{j \neq i} \frac{1}{\lambda_i - \lambda_j}
.
\]
\end{lemma}

\proofstart
This follows from the theory of partial fractions.

Alternatively, combine terms on the right over the common denominator
$\prod_i (x-\lambda_i)$. The numerator is
\[
\sum_i \prod_{j \neq i} \frac{x-\lambda_j}{\lambda_i-\lambda_j}
.
\]
This is a polynomial of degree $n-1$ which takes the value $1$ for
$x=\lambda_1,\ldots,\lambda_n$.
These $n$ values of $x$ are distinct (thinking of the $\lambda_i$'s as indeterminates),
so the numerator is identically $1$.
\proofend

{\bf Proof of Theorem \ref{one}.}
For general $q$, $s \in (\ZZ_q)^m$ put
\[
H_{q,s}(x,y) = \sum_{k \in \ZZ_q} \prod_i \frac{y-\omega_q^{ks_i}}{x-\omega_q^{k s_i}}
\]
so that
\[
\hodge_{L(q,s)}(x,y) = \frac{1}{q}\frac{y^n}{x^n} H_{q,s}(-1/y,1/x)
.
\]
Separate the sum for $H_{q,s}$
into pieces according to $\gcd(k,q)$ by putting
\[
H^\star_{d,s} 
=
\sum_{k \in \ZZstar_d} 
\prod_i \frac{y-\omega_d^{ks_i}}{x-\omega_d^{k s_i}}
\]
so that
\[
H_{q,s} = \sum_{d \divides q} H^\star_{d,s}
.
\]

To prove the theorem, we must show that if $a$ is hereditarily good mod
$r$ then
\[
H_{r^2t,\dub{(rta+1)}}
=
H_{r^2t,\dub{(-rta+1)}}
.
\]
Our strategy will be to show that for all $d \divides r^2t$ we have
\[
H^\star_{d,\dub{(rta+1)}}
=
H^\star_{d,\dub{(-rta+1)}}
.
\]

We dispose first of the case where $\dub{(rta+1)}$ is not univalent mod $d$.
This is taken care of by the assumption that $a$ is hereditarily good,
but things are not quite as straight-forward as you might be expecting,
because that condition deals with divisors of $r$, and here $d$ is any
divisor of $r^2t$.

We pass over the trivial cases $d=1,2$.
Mod any $d>2$, there is no overlap between $rta+1$ and $-(rta+1)$, so 
if $\dub{(rta+1)}$ is not univalent mod $d$ then neither is $rta$.

We pause for a lemma.

\begin{lemma}
For $d,\alpha,\beta \in \ZZ$, suppose $d \divides \alpha \beta$.
Let $d'=d/\gcd(d,\beta)$.
Then $d'  \divides  \alpha$ and 
\[
\forall \gamma \in \ZZ \; (d \divides \beta \gamma \iff d'  \divides  \gamma)
.
\]
\end{lemma}

\proofstart
Let $e=\gcd(d,\beta)$, so that $d'=d/e$.
\[
d \divides \alpha \beta \implies d'=d/e \divides \alpha \beta /e
\]
and $\gcd(d',\beta/e)=1$ so $d' \divides  \alpha$, and for any $\gamma$
\[
d  \divides  \beta \gamma \iff d'  \divides \beta/e \gamma \iff d'  \divides  \gamma.
\]
(Pretty standard stuff, admittedly.)
\proofend

So suppose $d \divides r^2t$ and $d \divides rt(a_i-a_j)$ for $i \neq j$.
Putting $\alpha=r$, $\beta=rt$, $\gamma=a_i-a_j$ in the lemma
we get
\[
d'=d/\gcd(d,rt)  \divides  \alpha =r
\]
and
\[
d'  \divides  \gamma = a_i - a_j
.
\]
This tells us that $a$ is not univalent mod $d'$,
but since by assumption it is good mod any divisor of $r$,
it must be reversible mod $d'$:
\[
a \eqmd{d'} -a
.
\]
By the lemma, this is equivalent to
\[
rta \eqmd{d} -rta
\]
hence
\[
rta+1 \eqmd{d} -rta+1
.
\]
From this we get
\[
H^\star_{d,\dub{(ra+1)}}
=
H^\star_{d,\dub{(-ra+1)}}
.
\]

So from here on we may assume that $\dub{(rta+1)}$ (and hence also
$\dub{(-rta+1)}$) is univalent mod $d$, with $d \divides r^2t$.

Returning for a moment to the case of $H^\star_{d,s}$ for general $d,s$,
suppose $s \in (\ZZstar_d)^m$
with all the $s_i$'s distinct mod $q$,
so that the mod-$d$ quotient $s_j \quot{d} s_i$ is defined for all $i,j$,
and different from $1$ for $i \neq j$.
With this restriction, for $k \in \ZZstar_d$ we have
\[
\prod_i \frac{y-\omega_d^{ks_i}}{x-\omega_d^{k s_i}}
=
\sum_i
\frac{y-\omega_d^{k s_i}}{x-\omega_d^{k s_i}}
\prod_{j \neq i}
\frac{y-\omega_d^{k s_j}}{\omega_d^{k s_i}-\omega_d^{k s_j}}
.
\]
So
\begin{eqnarray*}
H^\star_{d,s} 
&=&
\sum_{k \in \ZZstar_d}
\sum_i
\frac{y-\omega_d^{k s_i}}{x-\omega_d^{k s_i}}
\prod_{j \neq i}
\frac{y-\omega_d^{k s_j}}{\omega_d^{k s_i}-\omega_d^{k s_j}}
\\&=&
\sum_{l \in \ZZstar_d}
\frac{y-\omega_d^{l}}{x-\omega_d^{l}}
\sum_i
\prod_{j \neq i}
\frac{y-\omega_d^{l s_j/s_i}}{\omega_d^{l}-\omega_d^{l s_j \quot{d} s_i}}
\\&=&
\sum_{l \in \ZZstar_d}
Y_{d,s}(x,y,\omega_d^l)
,
\end{eqnarray*}
where
\[
Y_{d,s}(x,y,w)
=
\frac{y-w}{x-w}
\sum_i
\prod_{j \neq i}
\frac{y-w^{s_j \quot{d} s_i}}{w - w^{s_j \quot{d} s_i}}
.
\]

Now recall our notation
\[
\dub{s} = (s_1,-s_1,\ldots,s_m,-s_m)
.
\]
Assuming the entries of $\dub{s}$ are all invertible and distinct mod $d$,
\[
Y_{d,\dub{s}} =
\frac{(y-w)(y-w^{-1})}{(x-w)(x-w^{-1})}
\sum_i
\prod_{j \neq i}
\frac
{(y-w^{s_j \quot{d} s_i})(y-w^{-s_j \quot{d} s_i})}
{(w - w^{s_j \quot{d} s_i})(w - w^{-s_j \quot{d} s_i})}
.
\]

Specializing finally to the case at hand,
take
\[
s = rta+1 = (rt a_1+1,\ldots,rt a_m+1)
\]
so that
\[
\dub{s} = 
(rta_1+1,-rta_1-1,\ldots,rta_m+1,-rta_m-1)
,
\]
and assume that these entries are all distinct mod $d$.

Here comes the magic:
For any $d  \divides  r^2t$ we have
\[
s_j \quot{d} s_i \equiv_d  rt (a_j-a_i)+1
.
\]
As the entries of $\dub{s}$ are distinct mod $d$, putting
\[
x_i = w^{rta_i}
\]
we have
\[
w^{s_j \quot{d} s_i} = 
\frac{x_j}{x_i}w
\]
and
\[
w^{-s_j \quot{d} s_i} = \frac{x_i}{x_j}w^{-1}
\]
so
\[
Y_{d,\dub{(ra+1)}}
=
\frac{(y-w)(y-w^{-1})}{(x-w)(x-w^{-1})}
\sum_i
\prod_{j \neq i}
\frac
{(y-\frac{x_j}{x_i}w)(y-\frac{x_i}{x_j}w^{-1})}
{(w-\frac{x_j}{x_i}w)(w-\frac{x_i}{x_j}w^{-1})}
.
\]
Setting $u=y/w$, $v=w^{-2}$, we get
\[
Y_{r^2t,\dub{(rta+1)}} =
\frac{(y-w)(y-w^{-1})}{(x-w)(x-w^{-1})}
F((x_1,\ldots,x_m),u,v)
\]
where
\[
F((x_1,\ldots,x_m),u,v)
=
\sum_{i=1}^m \prod_{j \neq i}
\frac
{(u-\frac{x_j}{x_i})(u-\frac{x_i}{x_j}v)}
{(1-\frac{x_j}{x_i})(1-\frac{x_i}{x_j}v)}
.
\]

Simultaneously we have
\[
Y_{r^2t,\dub{(-rta+1)}}
=
\frac{(y-w)(y-w^{-1})}{(x-w)(x-w^{-1})}
F((1/x_1,\ldots,1/x_m),u,v)
.
\]
In the next section we will prove the identity
\[
F((x_1,\ldots,x_m),u,v)
=
F((1/x_1,\ldots,1/x_m),u,v)
,
\]
from which we conclude
\[
H^\star_{d,\dub{(rta+1)}}
=
H^\star_{d,\dub{(-rta+1)}}
.
\]

We have now established this last equality for every $d \divides r^2t$, so
\[
H_{d,\dub{(ra+1)}}
=
H_{d,\dub{(-ra+1)}}
.
\mathproofend
\]

\section{The main identity}

Define the rational function
\begin{eqnarray*}
F((x_1,\ldots,x_m),u,v)
&=&
\sum_{i=1}^m \prod_{j \neq i}
\frac
{(u-\frac{x_j}{x_i})(u-\frac{x_i}{x_j}v)}
{(1-\frac{x_j}{x_i})(1-\frac{x_i}{x_j}v)}
\\&=&
\sum_{i=1}^m \prod_{j \neq i}
\frac
{(x_i u-x_j)(x_j u-x_i v)}
{(x_i  -x_j)(x_j  -x_i v)}
.
\end{eqnarray*}
Now look at what you get by replacing the variables
$x_1,\ldots,x_m$ by their reciprocals:
\begin{eqnarray*}
G((x_1,\ldots,x_m),u,v)
&=&
F((1/x_1,\ldots,1/x_m),u,v)
\\&=&
\sum_{i=1}^m \prod_{j \neq i}
\frac
{(u-\frac{x_i}{x_j})(u-\frac{x_j}{x_i}v)}
{(1-\frac{x_i}{x_j})(1-\frac{x_j}{x_i}v)}
\\&=&
\sum_{i=1}^m \prod_{j \neq i}
\frac
{(x_j u-x_i)(x_i u-x_j v)}
{(x_j  -x_i)(x_i  -x_j v)}
\end{eqnarray*}

\begin{prop}
$F=G$.
\end{prop}

\proofstart
The right way to prove this identity is presumably via invariant theory.
(Or maybe it's just somehow obvious?)
But here we are going to prove it by considering the two sides
as rational functions of $v$,
expanding their individual terms in partial fractions,
and seeing that the parts on the two sides agree.

The tricky case turns out to be the polynomial term,
corresponding to the pole at $v=\infty$.
We'll deal with that later, after we address the finite poles.

Let's look at the case $m=3$, which is sufficient to show
what is going on.
Each side of the identity has three terms.
On the left the first term is
\[
\frac{(x_1 u - x_2)(x_2 u - x_1 v)(x_1 u - x_3)(x_3 u - x_1 v)}
{(x_1 - x_2)(x_2 - x_1 v)(x_1 - x_3)(x_3 - x_1 v)}
.
\]
This term is the only one on the left with non-zero residue
at $v=x_2/x_1$, and its residue there is
\begin{eqnarray*}
&&
\frac{(x_1 u - x_2)(x_2 u - x_1 x_2/x_1)(x_1 u - x_3)(x_3 u - x_1 x_2/x_1)}
{(x_1 - x_2)(-x_1)(x_1 - x_3)(x_3 - x_1 x_2/x_1)}
\\&=&
\frac{(x_1 u - x_2)x_2(u - 1)(x_1 u - x_3)(x_3 u - x_2)}
{(x_1 - x_2)(-x_1)(x_1 - x_3)(x_3 - x_2)}
.
\end{eqnarray*}
On the right the only term with a non-zero residue at $v=x_2/x_1$ is the
second term, namely
\[
\frac
{(x_1 u - x_2)(x_2 u - x_1 v)(x_3 u - x_2)(x_2 u - x_3 v)}
{(x_1 - x_2)(x_2 - x_1 v)(x_3 - x_2)(x_2 - x_3 v)}
,
\]
and the residue there is
\begin{eqnarray*}
&&
\frac{(x_1 u - x_2)(x_2 u - x_1 x_2/x_1)(x_3 u - x_2)(x_2 u - x_3 x_2/x_1)}
{(x_1 - x_2)(-x_1)(x_3 - x_2)(x_2 - x_3 x_2/x_1)}
\\&=&
\frac{(x_1 u - x_2)x_2(u - 1)(x_3 u - x_2)(x_1 u - x_3 )}
{(x_1 - x_2)(-x_1)(x_3 - x_2)(x_1 - x_3)}
,
\end{eqnarray*}
which is the same as we found for the left side.

In this way we see that the residues of $v$ at the finite poles all match
between left and right.
That leaves the pole at $v=\infty$.
Taking the limit $v \goesto \infty$ of
\[
F((x_1,\ldots,x_m),u,v)
=
\sum_{i=1}^m \prod_{j \neq i}
\frac
{(x_i u-x_j)(x_j u-x_i v)}
{(x_i  -x_j)(x_j  -x_i v)}
\]
yields
\[
\sum_{i=1}^m \prod_{j \neq i}
\frac
{(x_i u-x_j)}
{(x_i  -x_j)}
.
\]
In the next section, we will prove that this limit is $1+u+\ldots+u^{m-1}$,
which as it is independent of
$(x_1,x_2,x_3)$ must agree with the limit of
\[
G((x_1,\ldots,x_m),u,v)=F((1/x_1,\ldots,1/x_m),u,v)
,
\]
so the residues at $v=\infty$ of the two sides of your identity match,
and the proof is complete.
Well, it's not a proof, exactly, since it doesn't really explain what is
going on there.
Call it a `verification', which persuades us that the identity is true,
at least when coupled with a symbolic computation checking the identity
up through $m=4$.
\proofend

\section{The subsidiary identity}

Define the rational function
\[
f((x_1,\ldots,x_m),u)
=
\sum_{i=1}^m \prod_{j \neq i} \frac{u-\frac{x_j}{x_i}}{1-\frac{x_j}{x_i}}
=
\sum_{i=1}^m \prod_{j \neq i} \frac{x_i u-x_j}{x_i-x_j}
.
\]

\begin{prop}
\[
	f((x_1,\ldots,x_m),u) = 1+u+\ldots+u^{m-1}
	.
\]
\end{prop}

\proofstart
We use induction on $m$.
The cases $m=0,1$ are trivial,
and $m=2$ is so easy as not to illustrate the method.
So we will look at the case $m=3$, and take that as representative.
We want to show that
\[
	\frac{(x_1 u - x_2)(x_1 u - x_3)}{(x_1-x_2)(x_1-x_3)}
	+
	\frac{(x_2 u - x_1)(x_2 u - x_3)}{(x_2-x_1)(x_2-x_3)}
	+
	\frac{(x_3 u - x_1)(x_3 u - x_2)}{(x_3-x_1)(x_3-x_2)}
	= 1 + u + u^2
	.
\]
Expand the terms on the left in partial fractions with respect to the variable $x_3$.
The possible poles are at $x_3=x_1$, $x_3=x_2$, and $x_3 = \infty$.
For the coefficient of $\frac{1}{x_3-x_1}$ we get
\[
	-\frac{(x_1 u - x_2)(x_1 u - x_1)}{x_1-x_2} + 0 + 
	\frac{(x_1 u - x_1)(x_1 u - x_2)}{x_1 - x_2}
	.
\]
So this coefficient vanishes (as it would have to, if our identity is to hold).
Similarly for the coefficient of $\frac{1}{x_3-x_2}$.
This leaves the pole at $x_3=\infty$.
Taking the limit of the terms on the
left as $x_3 \goesto \infty$, we get
\[
	\frac{x_1 u -x_2}{x_1-x_2} + \frac{x_2 u - x_1}{x_2 - x_1} + u^2
	= f((x_1,x_2),u) + u^2 = 1+u+u^2
	,
\]
where in the last step we are using the induction hypothesis.
\proofend

\section{Open questions}

\begin{enumerate}
\item
What is the right way to prove these two identities?
\item
LMR showed that in their construction, the full Hodge series agree
just if they agree after setting $w=0$.
Surely we can prove this algebraically.
\item
Is the condition in Theorem \ref{one} for Hodge-equivalence of LMR groups
necessary as well as sufficient?
If true, this might not be so hard to prove.
To start with, we could prove that what works for $t=1$ works for any $t$.
\item
It seems that the representation-theoretic proof of LMR might give
an explicit matchup between spaces of invariant forms.
Can we extract such a matchup from the algebra in the proof
of Theorem \ref{one}?
\item
Not all Hodge-isospectral pairs emerge directly from the LMR construction.
For example, you append a $0$ to the list $rta+1$ on both sides,
or put in everything congruent to $2$ mod $rt$.
It's tempting to figure out just what variations are possible.
And then we could ask whether all possible pairs arise as variations of
this kind.
\item
Doyle and Rossetti
\cite{doylerossetti:orbv1}
conjectured that in spherical geometry or hyperbolic geometry,
spaces that are $p$-isospectral for all $p$ are almost
conjugate,
and hence isospectral for all natural operators.
The LMR examples show that this is false in spherical geometry,
but the hyperbolic case remains open,
and the intuition for this conjecture,
born in the hyperbolic case and incautiously extended to the spherical
case, remains more or less intact.
\end{enumerate}

\section*{Acknowledgements}
We thank Emilio Lauret for extensive correspondence.

\bibliography{hodge}

\def\cprime{$'$}
\begin{thebibliography}{10}

\bibitem{crass:sextic}
Scott Crass.
\newblock Solving the sextic by iteration: {A} study in complex geometry and
  dynamics.
\newblock {\em Experimental Mathematics}, 8, 1999, arXiv:math/9903111
  [math.DS].
\newblock \url{http://arxiv.org/abs/math/9903111}.

\bibitem{doylemcmullen:icos}
Peter Doyle and Curt McMullen.
\newblock Solving the quintic by iteration.
\newblock {\em Acta Mathematica}, 163:151--180, 1989.
\newblock
  \url{http://www.math.harvard.edu/~ctm/papers/home/text/papers/icos/icos.pdf}.

\bibitem{doylerossetti:orbv1}
Peter~G. Doyle and Juan~Pablo Rossetti.
\newblock {Laplace-isospectral hyperbolic 2-orbifolds are
  representation-equivalent (Version 1)}, 2011, arXiv:1103.4372v1 [math.DG].
\newblock \url{http://arxiv.org/abs/1103.4372v1}.

\bibitem{gilkey:metacyclic}
Peter~B. Gilkey.
\newblock On spherical space forms with metacyclic fundamental groups which are
  isospectral but not equivariant cobordant.
\newblock {\em Compositio Math.}, 56:171--200, 1985.

\bibitem{vihart:pi}
Vi~Hart.
\newblock Pi is (still) wrong.
\newblock \url{https://www.youtube.com/watch?v=jG7vhMMXagQ}.

\bibitem{ikeda:subtle}
Akira Ikeda.
\newblock Riemannian manifolds $p$-isospectral but not $(p+1)$-isospectral.
\newblock In {\em Geometry of Manifolds}, pages 383--417. Academic Press, 1989.

\bibitem{klein:ikos}
Felix Klein.
\newblock {\em Vorlesungen {\"u}ber das {I}kosaeder und die {A}ufl{\"o}sung der
  {G}leichungen vom f{\"u}nften {G}rade}.
\newblock B.G. Teubner, 1884.
\newblock \url{http://www.math.dartmouth.edu/~doyle/docs/ikos/scan/ikos.pdf}.

\bibitem{lmr:lens}
Emilio~A. Lauret, Roberto~J. Miatello, and Juan~Pablo Rossetti.
\newblock Lens spaces isospectral on $p$-forms for every $p$, 2013,
  arXiv:1311.7167v2 [math.DG].
\newblock \url{http://arxiv.org/abs/1311.7167v2}.

\bibitem{molien}
Theodor Molien.
\newblock {\"{U}ber} die {Invarianten} der linearen {Substitutionsgruppen}.
\newblock {\em Sitzungber. K\"{o}nig. Preuss. Akad. Wiss. (J. Berl. Ber.)},
  52:1152--1156, 1897.
\newblock \url{http://books.google.com/books?id=EIxK-opAmJYC&pg=PA1152}.

\bibitem{pesce:strong}
Hubert Pesce.
\newblock Vari\'et\'es hyperboliques et elliptiques fortement isospectrales.
\newblock {\em J. Funct. Anal.}, 134(2):363--391, 1995.

\bibitem{stanley:molien}
Richard~P. Stanley.
\newblock Invariants of finite groups and their applications to combinatorics.
\newblock {\em Bull. Amer. Math. Soc. (N.S.)}, 1(3):475--511, 1979.
\newblock \url{http://math.mit.edu/~rstan/pubs/pubfiles/38.pdf}.

\end{thebibliography}
\bibliographystyle{hplain}
\end{document}